\documentclass[12pt,leqno,twocolumn]{article}
\usepackage{latexsym}
\usepackage{epsfig}

\makeatletter
\@addtoreset{equation}{section}

\makeatother

\begin{document}
\title{Sphere Packings in $3$ Dimensions}
\author{Thomas C.
Hales}

\maketitle

\newtheorem{remark}{Remark}[section]
\newtheorem{lem}{Lemma}[section] \newtheorem{thm}{Theorem}[section]
\newtheorem{cor}{Corollary}[section]
\newtheorem{pro}{Proposition}[section]
\newtheorem{defn}{Definition}[section]

The Kepler conjecture asserts that no packing of equal balls in three
dimensions can have density exceeding that of the face-centered cubic
packing.  The density of this packing is $\pi/\sqrt{18}$, or about 0.74048.

This conjecture was proved in 1998 [H].  In superficial terms, the
proof takes over 270 pages of mathematical text, extensive computer
resources, including 3 gigabytes of data, well over 40 thousand lines
of code, about $10^{5}$ linear programs each involving perhaps 200
variables and 1000 constraints.  The project spanned nearly a decade
of research, including substantial contributions from S.  Ferguson
[F].  This is not a {\it proof from the book}, in the sense of Erd\"os.

This lecture\footnote{This article is based on an Arbeitstagung
lecture in Bonn on June 14, 2001.} proposes a second-generation
proof of the Kepler conjecture, and describes the current status
of that second-generation proof.  Much of the inspiration for the
revised proof comes from the subject of {\it generative
programming}, and the lecture gives a brief description of that
topic as well.

There are various motivations for a second-generation proof.  One
is that a simple is more likely to be used as a model for other
proofs in discrete geometry.  One might try to build a proof of
the kissing number problem in 4 dimensions, or a proof of the
Kelvin problem (of finding the optimal partition of Euclidean
3-space into equal volumes) by similar methods.  Such an approach
is appealing only if the methods can be further simplified.

This lecture is primarily concerned with packings in three dimensions.
But it should be pointed out that excellent progress in higher
dimensions has been made in recent work of H.  Cohn and N.  Elkies.
They have produced the best known bounds on densities of packings in dimensions
4-36, and they conjecture that their methods will lead to sharp bounds
in dimensions 8 and 24 (corresponding to the $E_{8}$ and Leech
lattices respectively) [CE].

Two common decompositions of space are the Voronoi decomposition and
the Delaunay decomposition.  These two decompositions lead to bounds
of 0.755 and 0.740873, respectively (the second bound is not rigorous,
being based merely on numerical evidence) [M].  In the 1998 proof of
the Kepler conjecture, Euclidean space was partitioned according to a
hybrid of these two decompositions, and this led to the exact bound of
$\pi/\sqrt{18}$.

In the second-generation proof, a function is associated with each
center $\lambda$ of each ball.  Assume that all balls in the packing
have radius 1.  The functions $f_{\lambda}$ at $\lambda$ is the sum of
the volume of the Voronoi cell around $\lambda$ and a correction term
$\epsilon_{\lambda}$.  The correction term has the form
$$\epsilon_{\lambda} = \sum\limits_{T} \ \delta(T) + M \
\sum\limits_{S} \ \mu_{S,\lambda} \ (\sqrt{8} - w(S)),$$
where $T$ runs over all triangles with vertices at centers of balls,
and one vertex at $\lambda$, such that the edges have lengths at most
some fixed length $r \in [2, \sqrt{8}]$.  By definition,
$$\delta(T) = 2L(c) - L(a) - L(b)$$
where $L$ is a quadratic polynomial correction term.  Also, $a, b$ and
$c$ are the lengths of the edges of $T$, and $c$ is the edge opposite $\lambda$.

The term proportional to $M$ is a bit more difficult to describe.
$M$ is a constant.  The sum is over certain triangles $S$ with
vertices at the centers of balls with one vertex at $\lambda$.
(We omit the definition of these triangles.) Each of these
triangles $S$ has a distinguished edge of length $w(S)$.  The
multiplicity $\mu$ is $\mu_{S,\lambda} = -1$ if an endpoint of the
distinguished edge of $S$ is $\lambda$ and $\mu_{S,\lambda} = 2$
otherwise.

\noindent
{\bf Conjecture 1.} For a suitable quadratic polynomial $L$, constant
$M$, and parameter $r$, we have $f_{\lambda} \geq \nu_{0}$, where
$\nu_{0}$ is the volume of the rhombic dodecahedron (the Voronoi cell
of the face-centered cubic packing).

This conjecture implies the Kepler conjecture, because the correction terms
$\epsilon_{\lambda}$ cancel on average, when we sum over all centers of balls in
a packing.  The negative terms are designed to precisely cancel the positive
terms.  Thus, Conjecture 1 implies that Voronoi cells are, on average in some
large region, at least as large as those of the face-centered cubic packing (up
to a small error term depending on the size of the region), and this implies
that packings in general can be no denser than the face-centered cubic packing.

According to a classical result of Tarski (ca. 1930), for any
given $L, M$, and $r$ (with algebraic constants and coefficients),
the truth of the conjecture can be determined by a finite
algorithm.  Collins has an algorithm giving an improvement over
Tarski's original algorithm, but even that  is doubly exponential
in the number of variables, and so it is not practical for a
problem in over 100 variables [CJ].

Instead, the methods that one would use to prove this conjecture  are the same
as those that would be used in the original proof of the Kepler conjecture
(linear programming, interval arithmetic, automatic generation of graphs).
However, a solution of the Kepler conjecture obtained through Conjecture 1
would probably be much simpler than the original proof.  At the very least, the
minimization problem of Conjecture 1 is much simpler than the one considered in
the 1998 proof.

About 200 pages of the 1998 proof can be summarized in one sentence.  Break the
configuration space into a bunch of cases, and then show that inequalities
proved by interval arithmetic (together with a little bit of geometry) can be
applied in such a way to show that these cases are worse than the best known
close-packings.

Given that 200 pages can be summarized so easily, we should try to structure
the repetitive pattern into an algorithm and automate the proof.  Multi-simplex
optimization (to be explained in the lecture), combining linear programming
and  interval arithmetic, gives a general inequality-proving algorithm that can
be applied to automate much of the proof in this way.  Multi-simplex
optimization occurs in the 1998 proof, but only in a peripheral way.

Additional simplifications of the proof would come from generative
programming techniques.  For example, generative program could be
used to develop a general method for automated inequality proving
based on interval arithmetic [CZ].

\bigskip

\centerline{\bf REFERENCES}

\bigskip

\noindent
[CE] H.  Cohn and N.  Elkies, {\it New upper bounds on sphere packings I},
 math.MG/0110009.
\bigskip

\noindent [CJ] B.F.  Caviness and J.R.  Johnson (ed.) {\it
Quantifier Elimination and Cylindrical Algebraic Decomposition},
Springer-Verlag, 1998.
\bigskip

\noindent [CZ] K.  Czarnecki and U.  Eisenecker, {\it Generative
Programming: Methods, Tools and Applications}, Addison-Wesley, 2000.

\bigskip

\noindent [F] S.  Ferguson, {\it Sphere Packings V}, thesis, Univ.
of Michigan, 1997, math.MG/9811077.

\bigskip

\noindent
[H] T.  Hales, {\it The Kepler conjecture}, math.MG/9811078.

\bigskip

\noindent
[M] T.  Hales, S.  McLaughlin, {\it A proof of the dodecahedral conjecture},
math.MG/9811079.

\end{document}